\documentclass{article}

\usepackage[english]{babel}
\usepackage{subfigure}
\usepackage[letterpaper,top=2cm,bottom=2cm,left=3cm,right=3cm,marginparwidth=1.75cm]{geometry}
\usepackage{array}
\usepackage{mlmath}
\usepackage{appendix}
\usepackage{amsmath}
\usepackage{siunitx}
\usepackage{float}
\usepackage[section]{placeins}
\usepackage[linesnumbered,ruled]{algorithm2e}
\usepackage{amssymb}
\usepackage{graphicx}
\usepackage[colorlinks=true, allcolors=blue]{hyperref}
\title{Input gradient annealing neural network for solving low-temperature Fokker-Planck equations}
\author{Liangkai Hang, Dan Hu, Zhi-Qin John Xu}

\begin{document}
\maketitle

\begin{abstract}
We present a novel yet simple deep learning approach, called input gradient annealing neural network (IGANN), for solving stationary Fokker-Planck equations. Traditional methods, such as finite difference and finite elements, suffer from the curse of dimensionality. Neural network based algorithms are meshless methods, which can avoid the curse of dimensionality. Moreover, at low temperature, when directly solving a stationary Fokker-Planck equation with more than two metastable states in the generalized potential landscape, the small eigenvalue introduces numerical difficulties due to a large condition number. To overcome these problems, we introduce the IGANN method, which uses a penalty of negative input gradient annealing during the training. We demonstrate that the IGANN method can effectively solve high-dimensional and low-temperature Fokker-Planck equations through our numerical experiments.
\end{abstract}

\section{Introduction}
Rare events, which refer to processes with timescales far exceeding those accessible by brute force simulations~\cite{bonella2012theory}, often arise from high free energy barriers between different meta-stable states, such as chemical reactions and nucleation during phase transitions—ubiquitous phenomena in physical, chemical, and biological sciences \cite{zhang2016recent}. Studying rare events constitutes a central aim to elucidating the transition mechanisms between different meta-stable states, including characterizing the transition paths and the transition rates. The dynamics of many systems conform to a stochastic differential equation (SDE) as follows:
\begin{equation}
 \mathrm{d} \bm
{x}_t=\bm{f}\left(\bm{x}_t\right) \mathrm{d}t + \sqrt{2 \varepsilon} \bm{\sigma}(\bm{x}_t) \mathrm{d} \bm{W}_t, \label{eq:1}
\end{equation}
where $\bm{f}: \mathbb{R}^d \rightarrow \mathbb{R}^d $ is the dynamical driving force, $\bm{\sigma}({\bm{x}})  \in \mathbb{R}^{n \times m}$ is the diffusion matrix, $\bm{W}_t \in \mathbb{R}^m$ is an $m$-dimensional standard Brownian motion, $\varepsilon$ typically represents a small quantity that specifies the noise strength, and this quantity is often directly proportional to the system's temperature \cite{kubo1966fluctuation}. In physical systems, the driving force $\bm{f}$ is a potential force $-\nabla U$ \cite{lelièvre_stoltz_2016}, whereas in many biological systems with energy inputs, the force is not solely a potential force \cite{wang2011quantifying, wang2008potential}. High energy barrier implies that noisy strength $\varepsilon$ is very small. In this scenario, we refer to such a system as a low-temperature system. Then the probability density $p(\bm{x}, t)$ satisfies the Fokker-Planck equation \cite{risken1996fokker, wang2008potential},
\begin{equation}
    \frac{\partial p(\bm{x}, t)}{\partial t} = \mathcal{L} p(\bm{x}, t):=   - \nabla \cdot ( \bm{f}(\bm{x})p(\bm{x, t}) - \varepsilon \nabla \cdot (\bm{D}(\bm{x})p(\bm{x}, t))), \label{eq:2}
\end{equation}
where $\mathcal{L}$ is the Fokker-Planck operator and the matrix $\bm{D}(\bm{x}) = \bm{\sigma}(\bm{x})\bm{\sigma}(\bm{x})^{T} \in \mathbb{R}^{n \times n}$. The invariant distribution $p(\bm{x})$, governed by the stationary Fokker-Planck equation, i.e.,
\begin{equation}
    - \nabla \cdot ( \bm{f}(\bm{x})p(\bm{x, t}) - \varepsilon \nabla \cdot (\bm{D}(\bm{x})p(\bm{x}, t)))=0, \label{eq:sfp}
\end{equation}
is vital for elucidating the transition mechanisms \cite{wang2008potential, li2013quantifying}. 

~\\\
\noindent Traditional grid-based numerical methods, including finite difference \cite{berezin1987conservative, sepehrian2015numerical} and finite element methods \cite{naprstek2014finite, galan2007stochastic}, can be used to solve the steady-state Fokker-Planck equation to obtain the invariant distribution. However, discretizing Eq. \eqref{eq:sfp} over a finite spatial domain incurs exponentially increasing computational cost as the system dimension increases, limiting applicability of these methods. Furthermore, for systems with multiple meta-stable states, low system temperatures pose numerical challenges due to very small non-zero eigenvalues of the Fokker-Planck operator $\mathcal{L}$. These small eigenvalues are related to the low transition rate between meta-stable states. The small eigenvalues lead to a slow convergence rate when we solve Eq. \eqref{eq:2} or a high condition number of the discretized linear system of Eq. \eqref{eq:sfp}.

~\\\
\noindent The Monte Carlo approaches \cite{kikuchi1991metropolis, chen2017beating, chen2018efficient} can be used to overcome the curse of dimensionality. However, there are difficulties for this type of approaches to collect samples near the transition state. For gradient systems, i.e., the driving force $\bm{f}(\vx)$ is a potential force and the matrix $\bm{D}$ is isotropic, one can sample the system by introducing external potentials, such as umbrella sampling \cite{torrie1977nonphysical} and metadynamics \cite{laio2008metadynamics}. The invariant distribution can be obtained from the sampled distribution and external potentials. However, similar approaches cannot be used to find the invariant distribution for non-gradient systems.

~\\
\noindent Efficiently solving high-dimensional Fokker-Planck equations at low temperature  thus remains a challenge, especially for non-gradient systems. Deep neural networks (DNNs) provide a promising avenue to solve high-dimensional differential equations \cite{dissanayake1994neural,liu2020multi, raissi2019physics, weinan2018deep, weinan2021algorithms, zang2020weak,mod2022zhang,zhang2022multi}. Prior works have incorporated physical insights into DNNs for specific PDE contexts, including physical-informed neural networks (PINNs) \cite{dissanayake1994neural,raissi2019physics} and variational formulations \cite{weinan2018deep,zang2020weak,ming2021deep}. Theoretical works have made progress in understanding the mechanisms underlying deep learning, such as frequency principle of low-frequency bias \cite{xu_training_2018,xu2019frequency,CSIAM-AM-2-484,zhang2021linear,xu2022overview,liu2020multi}. For Fokker-Planck equations, neural network based methods have been developed by incorporating the normalization condition of the probability density function into loss functions to avoid trivial solutions \cite{xu2020solving}; by minimizing the loss to obtain a decomposition of the force field \cite{lin2022computing}; or by adding the entropy production rate (EPR) into loss function \cite{zhao2023epr}.

~\\\
\noindent Let $V= -\varepsilon \log p(\vx, t)$ represent the generalized potential. Instead of solving the original Fokker-Planck equation, we solve the Fokker-Planck equation in the general potential form, which can be written as
\begin{equation}
    \nabla V(\vx)^{T}(\bm{f}(\vx) - \varepsilon \nabla \cdot \bm{D}(\vx) + \bm{D}(\vx)\nabla V(\vx)) - \varepsilon \nabla \cdot (\bm{f}(\vx) - \varepsilon \nabla \cdot \bm{D}(\vx) + \bm{D}(\vx) \nabla V(\vx)) = 0. \label{eq: fokker-planck potential}
\end{equation}
The reason for this will be explained in the Preliminary section.

~\\\
\noindent In numerical experiments of solving Eq. \eqref{eq: fokker-planck potential} by a PINN, we find that for a low-temperature Fokker-Planck equation, although the loss is very small, the learned solution still has a very large difference with the exact solution. The key reason is that as the temperature goes to zero, there exists a trivial solution satisfying $\nabla V(\vx) = 0$. In experiments, we find that there are many flat segments in the learned potential function for low temperature systems, i.e., $\nabla V(\vx) = 0$. To address this challenge during training, we introduce an Input Gradient Annealing Neural Network (IGANN), which is penalized by  $-\nabla V(\vx)$ in the loss function. This term can avoid the trivial solution with $\nabla V(\vx) = 0$. To ensure the learning can converge to the exact solution, this penalty term gradually anneals during the training. To show the effectiveness of IGANN, we conducted multiple experiments, including high-dimensional and low-temperature problems, such as an eight-dimensional system with low-temperature of $100 k_{B}T$ energy barrier.

~\\\
\noindent This paper is structured as follows. Section 2 will review the Fokker-Planck equation and discusses challenges in numerical simulation, along with related deep learning approaches. Section 3 then delineates our proposed method. Section 4 demonstrates and evaluates the method on gradient and non-gradient systems. Finally, Section 5 presents conclusions.

\section{Preliminary}
In this section, we present some preliminary knowledge about the Fokker-Planck equation and offer explanations for the issues that arise when using traditional numerical methods for solving it.

\subsection{Fokker-Planck equation in the generalized potential form}
Instead of solving the stochastic differential equation (SDE) \eqref{eq:1} to determinate a particle's trajectory, we focus on the average behavior of a statistical ensemble of Brownian particles, i.e., the probability density function $p(\mathbf{x, t})$, which satisfies the Fokker-Planck equation (FPE) \eqref{eq:2}. The invariant distribution of the steady-state, denoted as $p_{ss}(\vx)$, satisfies:
\begin{equation}
     \nabla \cdot \mathbf{J} = 0, \label{eq:3}
\end{equation}
where the probability flux $\mathbf{J}$ is defined as $\mathbf{J} = \bm{f}(\vx)p_{ss}(\mathbf{\vx}) - \varepsilon \nabla \cdot (\bm{D}(\vx) p_{ss}(\mathbf{\vx}))$.

Rare event dynamics becomes important only in the following case: the probabilities of visiting some stable states are not small, and the transition times to escape from some of these stable states to visit the other states are not too long. Studying rare event dynamics involves three basic goals.
\begin{figure}[htbp]
	\centering
	\includegraphics[width=0.4\linewidth]{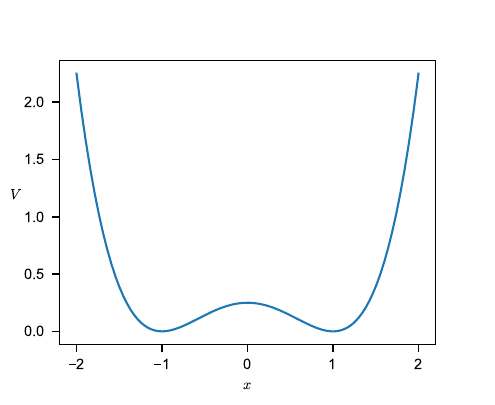}
	\caption{The generalized potential $V(x)$.}
	\label{fig: trueV}
\end{figure}

\begin{figure}[htbp]
	\centering
	\includegraphics[width=0.8\linewidth]{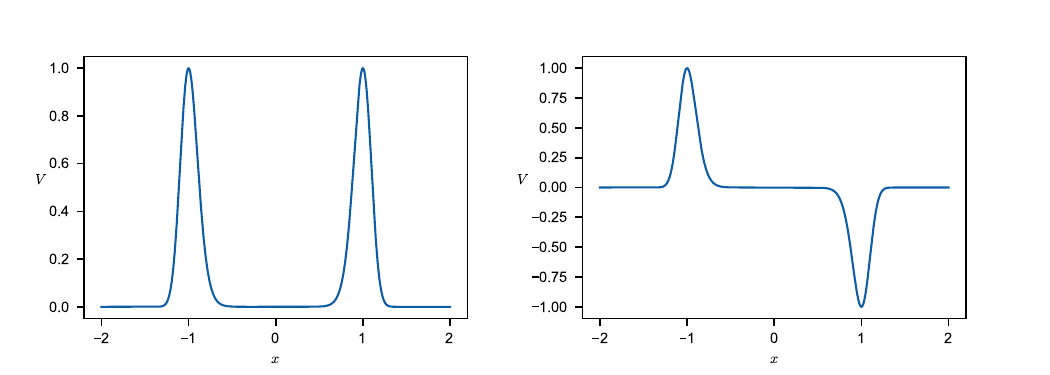}
	\caption{Eigenfunctions of the zero eigenvalue (left) and the first nonzero (right). Note that the left one is the steady-state solution of the Fokker-Planck equation up to a normalization constant.}
	\label{fig:eigenfunction}
\end{figure}

First, determine the probability of visiting each meta-stable state. This is simply the integral of the steady-state solution of the Fokker-Planck equation within each well. The peak values of the steady-state probability density function or the potentials within the wells, can roughly provide these probabilities. The first eigenvalue of the Fokker-Planck operator (denoted as $\mathcal{L}$), is zero and the corresponding eigenfunction is the steady-state solution. The second eigenvalue of $\mathcal{L}$, corresponding to the rare event dynamics, is very small (as discussed in the next subsection), and the associated eigenfunction has opposite signs in the example of two wells. As shown in Fig. \ref{fig:eigenfunction}, this system is a one-dimensional gradient system with a driving force $\bm{f} = -\nabla V(x)$. The generalized potential $V(x) = -\frac{x^{2}}{2} - \frac{x^{4}}{4}$, with an isotropic matrix $\bm{D} = \mathbf{I}$, and a noise strength $\epsilon=0.02$. The generalized potential $V$ is shown in Fig. \ref{fig: trueV}. Since the second eigenvalue is very small, it is hard to get rid of the second eigenfunction in the neural network solution of the Eq. \eqref{eq:3}, because the numerical error of the solution is usually large but the residual on the LHS of the Eq. \eqref{eq:3} is too small. Meanwhile, the contribution of the second eigenfunction to the neural network solution can significantly change the heights of the peaks in the probability density function, leading to a large error in estimating the probabilities of visiting the stable states.

Second, identify the transition mechanism(s) between meta-stable states. This is usually achieved by determining the transition path, which can be obtained by solving Eq. \eqref{eq:3} using string methods \cite{weinan2002string} and other numerical methods.

Third, estimate the transition rate. This is usually obtained either by evaluating the probability density at the saddle point (the transition point on the transition path) or by calculating the energy barrier (the energy difference between the steady-state and the saddle point). For a more accurate estimation, information from the Hessian matrix of the potential $V$ is needed.

Note that the probability density at the saddle point is usually very small (exponentially small). Thus the neural network solution of the Fokker-Planck equation is not able to characterize the smallness of the probability density function at the saddle point. This is why it is important to solve the steady-state Fokker-Planck equation using Eq. \eqref{eq: gpf}. let $p_{ss}(\vx) = \exp \left( \frac{-V(\vx)}{\varepsilon} \right)$, $\mathbf{l}_{ss}(\vx) = \bm{f}(\vx) - \varepsilon \nabla \cdot \bm{D}(\vx) + \bm{D}(\vx) \nabla V(\vx)$, then we have the equation for the generalized potential function $V(\vx)$:
\begin{equation}
    \varepsilon \nabla \cdot \mathbf{l}_{ss}(\vx) - \mathbf{l}_{ss}(\vx) \cdot \nabla V(\vx) = 0. \label{eq: gpf}
\end{equation}
Alternatively, we can decompose $\bm{f}(\vx)$ using the following formulation:
\begin{equation}
    \bm{f}(\vx) = -\bm{D}(\vx) \nabla V(\vx) + \varepsilon \nabla \cdot \bm{D}(\vx) + \mathbf{l}_{ss}(\vx),
\end{equation}
where $\mathbf{l}_{ss}(\vx)$ satisfies Eq. \eqref{eq: gpf}.
To illustrate the fundamental principles of the generalized potential landscape, consider a gradient system with
\begin{equation}
    \bm{f}(\vx) = -\nabla U(\vx), \, \, \bm{D}= \mathbf{I},
\end{equation}
it follows that $p_{ss}(\vx) = \frac{1}{Z}\exp \left(\frac{-U(\vx)}{\varepsilon} \right)$, implying that the generalized potential $V(\vx)$ is equivalent to the driving potential $U(\vx)$ up to a constant:
\begin{equation}
    V(\vx) = U(\vx) + \varepsilon \log Z.
\end{equation}
In biological systems, the generalized potential $V(\vx)$ is related to the global quasi-potential, which provides insights into the global behavior of the system \cite{zhou2016construction}.

\subsection{Limitations of traditional methods}
Grid-based methods such as finite difference (FDM) and finite element (FEM) face challenges of exponential scaling with dimensionality due to spatial discretization, which limits their applicability for high-dimensional problems. Moreover, a small $\varepsilon$ leads to additional challenges in traditional numerical methods, fundamentally stemming from a very small non-zero eigenvalue. 

If we assume that the Fokker-Planck operator possesses a discrete spectrum, then it is characterized by eigenfunctions $p_0, p_1, p_2, \ldots,$ corresponding to distinct eigenvalues $\lambda_0, \lambda_1, \lambda_2, \ldots$. Here, $\lambda_0$ is the zero eigenvalue and $p_0(\vx)$ is the steady-state solution. If the system satisfies the detailed balance condition, the Fokker-Planck operator $\mathcal{L}$ can be brought into a self-adjoint operator with real and non-negative eigenvalues. However, in most cases, $\mathcal{L}$ cannot be brought into a self-adjoint form. Under these circumstances, the first nonzero eigenvalue $\lambda_1$ is positive while others are complex-valued with a non-negative real part. \cite{matkowsky1981eigenvalues, holcman2014oscillatory, krein1948linear, blum1996variational} In the small noise limit, i.e., $\varepsilon \ll 1$, $\lambda_1$ is asymptotically related to the mean first passage time $\Bar{\tau}$ by the following equation:
\begin{equation}
  \lambda_1 \sim \frac{1}{\Bar{\tau}}.
\end{equation}

\begin{figure}[htbp]
	\centering
	\includegraphics[width=0.6\linewidth]{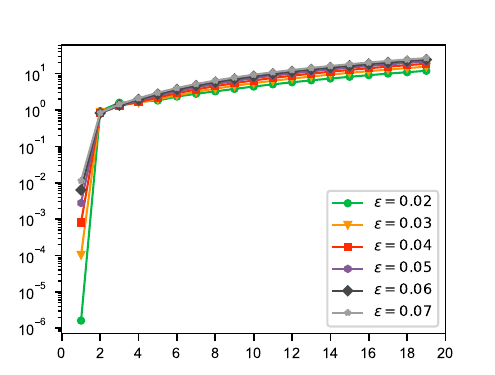}
	\caption{The magnitudes of the first twenty (in ascending order) non-zero eigenvalues of the matrix $\mathbf{A}$ are shown for varying epsilon values. The zero eigenvalue is not shown.}
	\label{eigenvalue}
\end{figure}

As the value of $\varepsilon$ approaches infinitesimally small magnitudes, $\lambda_1$ decays exponentially fast. For a visual illustration of how temperature affects this small eigenvalue, let us examine the one-dimensional system described previously. The computation domain is set to $[-2, 2]$. The eigenvalues of $\mathcal{L}$ can be computed within this domain using the finite difference method with zero-flux boundary conditions, i.e., $\mathbf{J}(\vx) \cdot \mathbf{n}=0$. Twenty eigenvalues are shown in Fig. \ref{eigenvalue}, where, as $\varepsilon$ gradually decreases, the second smallest eigenvalue $\lambda_1$ approaches zero. 

In this scenario, traditional methods face significant challenges. A prevalent strategy to obtain the invariant measure involves solving the steady-state Fokker-Planck equation. To achieve this, we discretize the operator of the Fokker-Planck equation using the finite difference method and specify a no-flux boundary condition. Subsequently, we obtain the following system of linear equations:
\begin{equation}
    \mathbf{A} \mathbf{p} = 0,
\end{equation}
accompanied by the  normalization condition:
\begin{equation}
   \int_{\Omega} p(\vx) \mathrm{d}\vx = 1,
\end{equation}
where $\Omega$ represents the computation domain.

To solve this system of linear equations, it is customary to set a particular element to 1, move the corresponding column to the right-hand side, and eliminate a row. This process results a new matrix $\widetilde{\mathbf{A}}$ and a new system of linear equations, given by $\widetilde{\mathbf{A}}\mathbf{p}=\mathbf{b}$. Due to the presence of the small eigenvalue $\lambda_1$ in the original matrix $\mathbf{A}$, the condition number of the new matrix $\widetilde{\mathbf{A}}$ becomes exceedingly high. This significantly increases the computational challenges in solving the new system of linear equations.

Monte Carlo approaches, enhanced by sampling methods such as umbrella sampling and metadynamics, are effective for gradient systems, but they are unsuitable for non-gradient systems. In practical scenarios, where we set the diffusion matrix $D(\vx)$ to be the identity matrix $I$ and the dimension of $\vx$ is very high, it becomes numerically difficult to perform the decomposition by solving the steady-state Fokker-Planck equation. Instead, one might be just interested in the probability density function $\pi(\mathbf{\xi})$, where $\mathbf{\xi}(\vx)$ is a lower-dimensional reaction coordinate. 
\begin{equation}
\pi(\mathbf{\xi})=\frac{1}{Z} \int_{\Omega} \exp \left(-\frac{V(\vx)}{\varepsilon}\right) \delta(\mathbf{\xi}(\vx)-\mathbf{\xi}) \mathrm{d}\vx.
\end{equation}
Since we do not have an explicit solution of $V(\vx)$ and the integral is in a high dimensional space, we wish to find $\pi(\mathbf{\xi})$ using a sampling method. In the subsequent discussion, we will briefly introduce the umbrella sampling method and discuss its inapplicability to non-gradient systems.

For a gradient system with $\bm{f}(\vx) = -\nabla U(\vx)$, we introduce an external potential $U^{b}(\vx)$. Consequently, the dynamical force becomes $-\nabla(U(\vx)+U^{b}(\vx))$. Thus, the biased density function $\pi^{b}(\mathbf{\xi})$ satisfies
\begin{equation}
\begin{aligned}
\pi^b(\mathbf{\xi}) & =\frac{1}{Z^b} \int \exp \left(-\frac{U(\vx)+U^b(\mathbf{\xi}(\vx))}{\varepsilon}\right) \delta(\mathbf{\xi}(\vx)-\mathbf{\xi}) \mathrm{d} \vx \\
& =\frac{1}{Z^b} \exp \left(-\frac{U^b(\mathbf{\xi})}{\varepsilon}\right) \int \exp \left(-\frac{U(\vx)}{\varepsilon}\right) \delta(\mathbf{\xi}(\vx)-\mathbf{\xi}) \mathrm{d} \vx \\
& =\frac{Z}{Z^b} \exp \left(-\frac{U^b(\mathbf{\xi})}{\varepsilon}\right) \pi(\mathbf{\xi}).
\end{aligned}
\end{equation}
Therefore, the free energy $\bar{u}(\mathbf{\xi}):=-\varepsilon \log \pi(\mathbf{\xi})$ can be calculated from the biased free energy $\bar{u}^{b}(\mathbf{\xi})$
\begin{equation}
    \bar{u}(\mathbf{\xi}) = \bar{u}^{b}(\mathbf{\xi}) - U^b(\mathbf{\xi}) + \varepsilon \log(\frac{Z}{Z^b}).
\end{equation}
However, when $\bm{f}(\vx)$ is not a gradient, namely, $\mathbf{l}(\vx) \neq 0$, the situation is different. Originally, we have
\begin{equation}
    \bm{f}(\vx) = -\nabla V(\vx) + \mathbf{l}(\vx),
\end{equation}
\begin{equation}
    \mathbf{l}(\vx) \cdot \nabla V(\vx) - \epsilon \nabla \cdot \mathbf{l}(\vx) = 0.
\end{equation}
When the external potential $U^b(\vx)$ is applied, the decomposition becomes
\begin{equation}
    \bm{f}(\vx) - \nabla U^b(\mathbf{\xi}(\vx)) = -\nabla V^b(\vx) + \mathbf{l}^b(\vx),
\end{equation}
\begin{equation}
    \mathbf{l}^b(\vx) \cdot \nabla V^b - \epsilon \nabla \cdot \mathbf{l}^b (\vx) = 0,
\end{equation}
obviously, we have 
\begin{equation}
    V(\vx) \neq V^{b}(\vx) - U^{b}(\mathbf{\xi}(\vx)):= R(\vx),
\end{equation}
denoting $T(\vx):= V(\vx) - R(\vx)$, it satisfies $\mathbf{l}(\vx) - \mathbf{l}^b(\vx) = \nabla T(\vx)$. In this case, we have
\begin{equation}
\begin{aligned}
\pi(\mathbf{\xi}) & =\frac{1}{Z} \int \exp \left(-\frac{V(\mathbf{x})}{\varepsilon}\right) \delta(\mathbf{\xi}(\mathbf{x})-\mathbf{\xi}) \mathrm{d} \vx \\
& =\frac{1}{Z} \exp \left(\frac{U^b(\mathbf{\xi})}{\varepsilon}\right) \int \exp \left(-\frac{V^b(\vx)+T(\vx)}{\varepsilon}\right) \delta(\mathbf{\xi}(\vx)-\mathbf{\xi}) \mathrm{d} \vx \\
& =\frac{Z^b}{Z} \exp \left(\frac{U^b(\mathbf{\xi})}{\varepsilon}\right)\left\langle\exp \left(-\frac{T(\vx)}{\varepsilon}\right) \delta(\mathbf{\xi}(\vx)-\mathbf{\xi})\right\rangle_b,
\end{aligned}
\end{equation}
where $\left\langle \cdot \right\rangle_b$ denotes the expectation under the biased probability density function. In principle, one can derive a nonlinear second-order partial differential equation (PDE) for $T(\vx)$. However, when the dimension of the dynamical system is high, it becomes impractical to determine $T(\vx)$ by solving the PDE of $T(\vx)$. Consequently, we must find an alternative method to estimate the expectation $\left\langle\exp \left(-\frac{T(\vx)}{\varepsilon}\right) \delta(\mathbf{\xi}(\vx)-\mathbf{\xi})\right\rangle_b$ using the sample data.

\section{Input gradient annealing network}
\begin{figure}[htbp]
	\centering
	\includegraphics[width=0.9\linewidth]{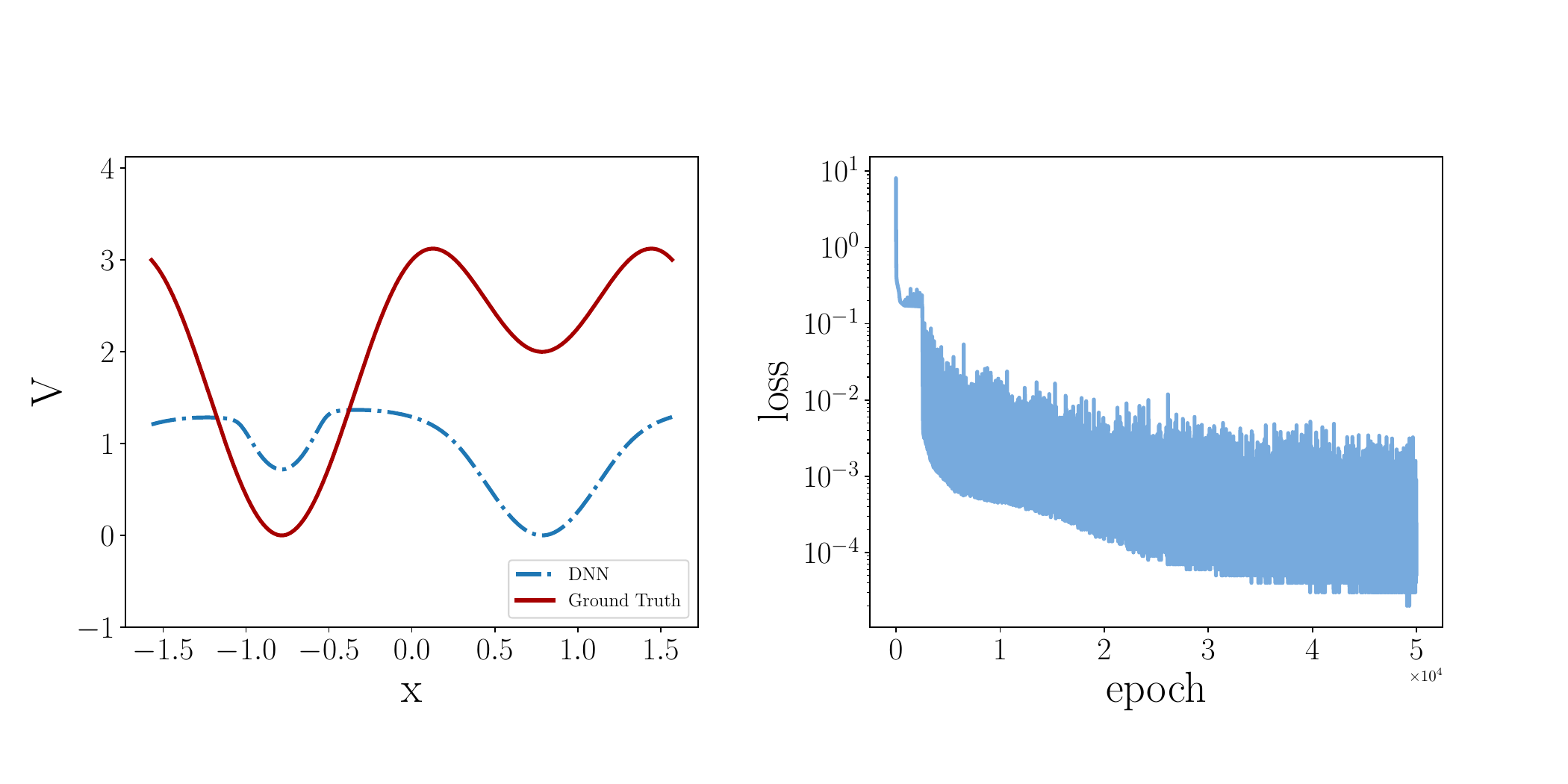}
	\caption{Minimizing the loss function \eqref{pinn_base}.}
	\label{1d_PINN_Fail}
\end{figure}

To effectively determine the generalized potential governed by the FPE, we employ deep neural networks (DNNs) to parameterize the generalized potential $V$. The least squared loss function is defined as 
\begin{equation}
    L = \int_{\mathbb{R}^d} \left| \nabla V_{\vtheta}(\vx)^T \mathbf{l}_{\vtheta}(\vx) - \varepsilon \nabla  \mathbf{l}_{\vtheta}(\vx) \right|^2 \mathrm{d} \mu(\vx), \label{pinn_base}
\end{equation}
where $\mathbf{l}_{\vtheta} (\vx) = \bm{f}(\vx) +  \bm{D}(\vx) \nabla V_{\vtheta}(\vx) - \varepsilon \nabla \cdot \bm{D}(\vx)$, and $\mu$ is the probability measure. However, solely minimizing loss defined in Eq. \eqref{pinn_base} often fails to recover the exact potential for low temperature system. As illustrated in Fig. \ref{1d_PINN_Fail} for a 1D system, despite the loss reaching 5e-5, the neural network poorly approximates the exact potential. Experiments reveal that the gradient of network with respect to the input $\vx$ tends to zero locally. Thus, unlike regular penalties, we incorporate a negative term to prevent vanishing gradients. Specifically, the loss is then given by the following equation:
\begin{equation}
   L = \int_{\mathbb{R}^d} \left| \nabla V_{\vtheta}(\vx)^T \mathbf{l}_{\vtheta}(\vx) - \varepsilon \nabla  \mathbf{l}_{\vtheta}(\vx) \right|^2 \mathrm{d} \mu(\vx) - \frac{\beta}{d} \int_{\mathbb{R}^d} \left| \nabla V_{\vtheta}(\vx)\right|^2 \mathrm{d} \mu(\vx).
\end{equation}
$\beta$ is a hyperparameter that decays with each training step, following the update rule $\beta
\leftarrow \beta(1-\beta_{\text{decay}})$, where both the initial value $\beta_0$ and the decay coefficient $\beta_{\text{decay}}$ are hyperparameters. The descending property of $\beta$ is to ensure the learning towards minimizing the original loss \eqref{pinn_base} in the final stage of the training. In all numerical examples, we employ samples from Latin hypercube sampling (LHS) \cite{latin} or from the simulation of SDE \eqref{eq:1} to discretize the integral form of the loss function, i.e.,

\begin{equation}
   L = \frac{1}{N} \sum_{i=1}^{N}  \left| \nabla V_{\vtheta}(\vx_i)^T \mathbf{l}_{\vtheta}(\vx_i) - \varepsilon \nabla  \mathbf{l}_{\vtheta}(\vx_i) \right|^2  -  \frac{1}{N} \frac{\beta}{d} \sum_{i=1}^{N} \left| \nabla V_{\vtheta}(\vx_i)\right|^2.  \label{our_loss}
\end{equation}

\section{Related works}
In this section, we review a related method from Ref. \cite{lin2022computing}, which will be used for comparison. This work modeled the potential as $V_{\theta}(\vx) = \tilde{V}_{\theta}(\vx) + \sum_{i=1}^d \rho_i(x_i-c_i)^2$, where $\tilde{V}_{\theta}(\vx)$ is a neural network, $\rho_i$ and $c_i$ are learnable parameters, and the residual component  $\mathbf{l}(\vx)$ via neural network $\mathbf{l}_{\vtheta}(\vx)$, giving $\bm{f}_{\theta}(\vx) = -\bm{D}\nabla V_{\vtheta}(\vx) + \mathbf{l}_{\theta}(\vx)$. Network weights were optimized by minimizing the loss function:
\begin{equation}
    L = L^{dyn} + \lambda L^{con},
\end{equation}
where 
\begin{equation}
\begin{aligned}
L^{d y n} & =\frac{1}{d} \int_{\mathbb{R}^d}\left|\bm{f}(\vx)-\bm{f}_\theta(\vx)\right|^2 d \mu(\vx), \\
L^{c o n} & =\int_{\mathbb{R}^d}\left|\nabla V_\theta(\vx)^T \mathbf{l}_\theta(\vx)-\varepsilon \nabla \cdot \mathbf{l}_\theta(\vx)\right|^2 d \mu(\vx),
\end{aligned}
\end{equation}
where $\mu(\vx)$ represents a probability measure, $L^{\text{dyn}}$ enforces $\bm{f}_{\vtheta} \approx \bm{f}$, $L^{\text{con}}$ imposes constraint \eqref{eq: fokker-planck potential} on $\bm{f}$ decomposition, and $\lambda$ balances the loss terms. The integrals above in the training process are approximated by finite sums using data points sampled from the uniform distribution on a bounded domain, or a mixture of the uniformly sampled data points and those sampled from the numerical simulation of the SDE \eqref{eq:1}. For convenience in comparing results, we will refer to this method as the LLR method.

\section{Results}
To demonstrate the efficacy of our approach, we apply it to gradient and non-gradient systems across various temperatures. For the non-gradient systems, we compare our method with the approach described in \cite{lin2022computing}. In particular, we apply our method to several examples presented in \cite{lin2022computing}, including the two metastable states, the biochemical oscillation network, and a ten-dimensional system. Our method demonstrates good performance in these examples, showing comparable results to those obtained in \cite{lin2022computing}. Furthermore, we construct another high energy barrier example to test the robustness of our method. In this challenging scenario, our method maintains its effectiveness, while the method from \cite{lin2022computing} encounters some difficulties in capturing the correct dynamics. In all examples, we utilize four-layer fully-connected neural networks equipped with the GELU activation function and He uniform initialization \cite{he2015delving}. During the training process, we perform LHS in each epoch or simulate SDE \eqref{eq:1} to obtain samples. We also use a cosine annealing learning rate scheduler. After training, the accuracy of the predicted potential is evaluated using the relative root mean square error (rRMSE) and relative mean absolute error (rMAE):
\begin{equation}
    rRMSE = \frac{\sqrt{\sum_{i=1}^{N} \frac{1}{N}(V_{\vtheta}(\vx_i) - V(\vx_i))^2}}{\sqrt{\sum_{i=1}^{N} V^2(\vx_i)}},
\end{equation}

\begin{equation}
    rMAE = \frac{\sum_{i=1}^{N} \frac{1}{N} \left| V_{\vtheta}(\vx_i) - V(\vx_i) \right|}{\sum_{i=1}^{N} \left| V(\vx_i) \right|},
\end{equation}
where $V_{\vtheta}$ is the learned generalized potential, and $V$ is the analytic solution or the numerical solution obtained by traditional methods. To facilitate comparison across different cross-sections, both solutions $V_{\vtheta}$ and $V$ are shifted so that their minimum values are zero.
\subsection{Gradient systems}
\begin{figure}[H]
	\centering
	\includegraphics[width=1\linewidth]{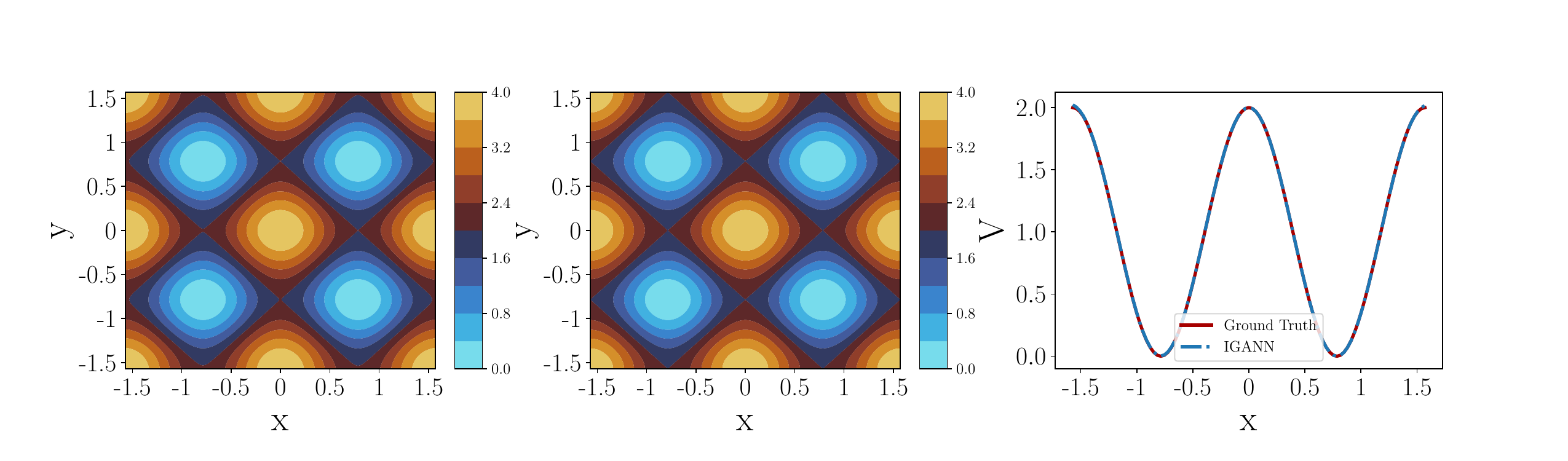}
	\caption{The learned potential $V_{\vtheta}$ compared with true potential for $\varepsilon = 0.04$ (the IGANN prediction (left), exact potential (middle) and central line $y=z=0$ (right).}
	\label{figure: ex4.1.3}
\end{figure}

We demonstrate the effectiveness of our method in a three-dimensional gradient system. The dynamical driving force is given by:
\begin{equation}
    \bm{f}: \vx = (x, y, z)^{T} \in \mathbb{R}^3 \rightarrow (4\sin(4x), 4\sin(4y), 4\sin(4z))^{T} \in \mathbb{R}^3,
\end{equation}
The corresponding dynamical system is:
\begin{equation}  
\left\{  
             \begin{array}{lr}  
             \mathrm{d} x = 4\sin(4x) \mathrm{d}t + \sqrt{2\varepsilon} \mathrm{d}W^x, \\  
             \mathrm{d} y = 4\sin(4y) \mathrm{d}t + \sqrt{2\varepsilon} \mathrm{d}W^y, \\
             \mathrm{d} z = 4\sin(4z) \mathrm{d}t + \sqrt{2\varepsilon} \mathrm{d}W^z, \\ 
             \end{array}  
\right.  
\end{equation}
where the computation domain is set to $[-\frac{\pi}{2}, \frac{\pi}{2}] \times [-\frac{\pi}{2}, \frac{\pi}{2}] \times [-\frac{\pi}{2}, \frac{\pi}{2}]$. For data collection, we gathered 2000 samples per epoch using LHS, and set $\beta=100$ with $\beta_{decay}=\num{3e-4}$. 
For $\varepsilon = 0.04$, considering the symmetry of the exact potential, we display only the cross-section at $z=0$ in Fig. \ref{figure: ex4.1.3}. The quantitative assessment, which includes results for other values of $\varepsilon$, is provided in Table \ref{table: ex4.1.3}. The results demonstrate good accuracy with very low error.
\FloatBarrier

\begin{table}[htbp]
\centering
\caption{Example \ref{sec: 4.1.3}}
\begin{tabular}{|p{3cm}<{\centering}|p{3cm}<{\centering}|p{3cm}<{\centering}|} \hline
&rRMSE (IGANN)&rMAE (IGANN)\\ \hline
$\varepsilon=0.1$ &  1.02\% & 0.518\% \\ \hline
$\varepsilon=0.07$&  0.667\%& 0.358\% \\ \hline
$\varepsilon=0.04$&  0.255\%& 0.154\%\\ \hline
\end{tabular}
\label{table: ex4.1.3}
\end{table}

\subsection{non-gradient systems}
In this section, we demonstrate the efficacy of our method in non-gradient systems.
\subsubsection{Two-dimensional systems} \label{sec: 4.2.1}

\begin{figure}[H]
	\centering
	\includegraphics[width=0.95\linewidth]{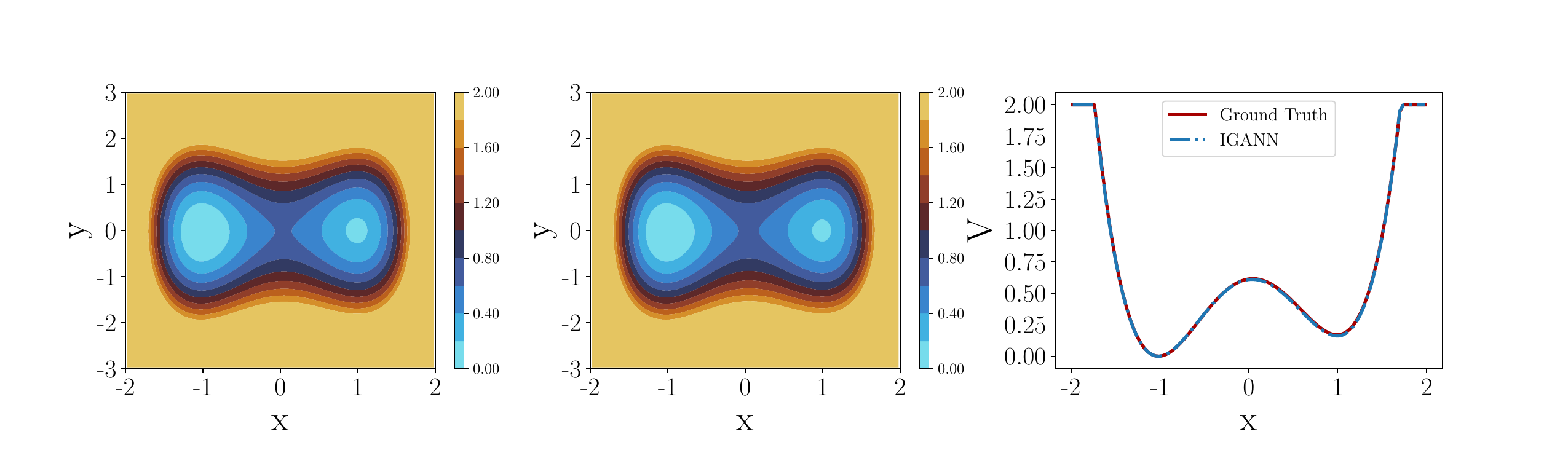}
	\caption{The learned potential $V_{\vtheta}$ compared with true potential for $\varepsilon = 0.1$ (the IGANN prediction (left), true potential (middle) and central line $y=0$ (right).}
	\label{fig: ex4.2.1 meta0.1}
\end{figure}

\begin{figure}[H]
	\centering
	\includegraphics[width=0.95\linewidth]{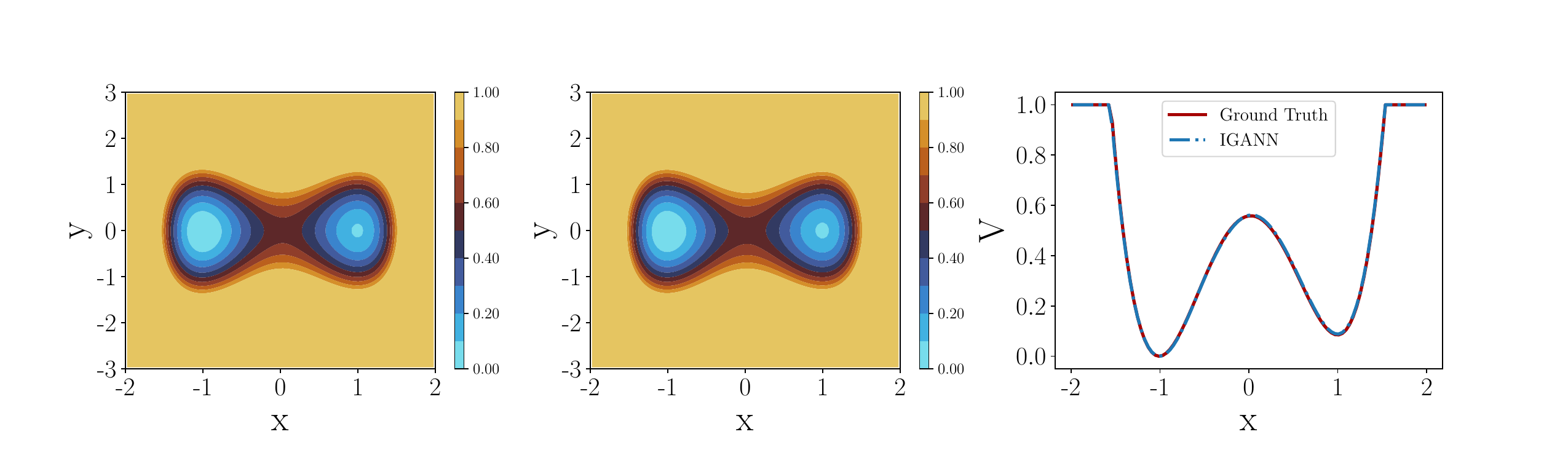}
	\caption{The learned potential $V_{\vtheta}$ compared with true potential for $\varepsilon = 0.05$ (the IGANN prediction (left), true potential (middle) and central line $y=0$ (right).}
	\label{fig: ex4.2.1 meta0.05}
\end{figure}

\textbf{Two metastable states dynamical system.}\label{ex:two metastable} We consider the following two-dimensional dynamical system described by 
\begin{equation}  
\left\{  
             \begin{array}{lr}  
             \mathrm{d} x = (\frac{1}{5} x(1-x^2) + y(1+\sin x)) \mathrm{d}t + \sqrt{\frac{\varepsilon}{5}} \mathrm{d}W^x, \\  
             \mathrm{d} y = ( -y + 2x(1-x^2)(1+\sin x)) \mathrm{d}t + \sqrt{2\varepsilon} \mathrm{d}W^y, \\
             \end{array}  
\right.  
\end{equation}
where the computation domain is set to $[-2, 2]\times[-3,3]$. We conduct computations for $\varepsilon=0.1, 0.05$ to show that our IGANN method can perform as well as the LLR method \cite{lin2022computing}. In this case, a small coefficient $\beta$ is sufficient, as larger $\beta$ values would increase the computational time. For $\varepsilon=0.1$, we set $\beta=5.$ with $\beta_{decay}=\num{4e-4}$. For $\varepsilon=0.05$ with $\beta_{decay}=\num{4e-4}$, we set $\beta=10$. The learned potentials using the IGANN method for $\varepsilon=0.1$ and $\varepsilon=0.05$ are illustrated in Fig. \ref{fig: ex4.2.1 meta0.1} and Fig. \ref{fig: ex4.2.1 meta0.05}, together with the finite difference solution. The cross section results and error comparisons are provided in Table \ref{table: ex4.2.1 meta}.
\FloatBarrier

\FloatBarrier

\begin{table}[htbp]
\centering
\caption{Example \ref{sec: 4.2.1} two metastable system}
\begin{tabular}{|p{2cm}<{\centering}|p{3cm}<{\centering}|p{3cm}<{\centering}|p{2.5cm}<{\centering}|p{2.5cm}<{\centering}|} \hline
 & $\varepsilon=0.1$ (IGANN)&  $\varepsilon=0.05$ (IGANN)& $\varepsilon=0.1$ (LLR)&  $\varepsilon=0.05$ (LLR) \\ \hline
rRMSE& 0.636\%& 0.66\%&  1.07\% $\pm$ 0.43\%& 1.94\% $\pm$  0.84\% \\ \hline
rMAE& 0.571\%& 0.605\%& 1.02\% $\pm$ 0.4\%& 1.93\% $\pm$ 0.9\%\\ \hline
\end{tabular}
\label{table: ex4.2.1 meta}
\end{table}

\begin{figure}[H]
	\centering
	\includegraphics[width=0.95\linewidth]{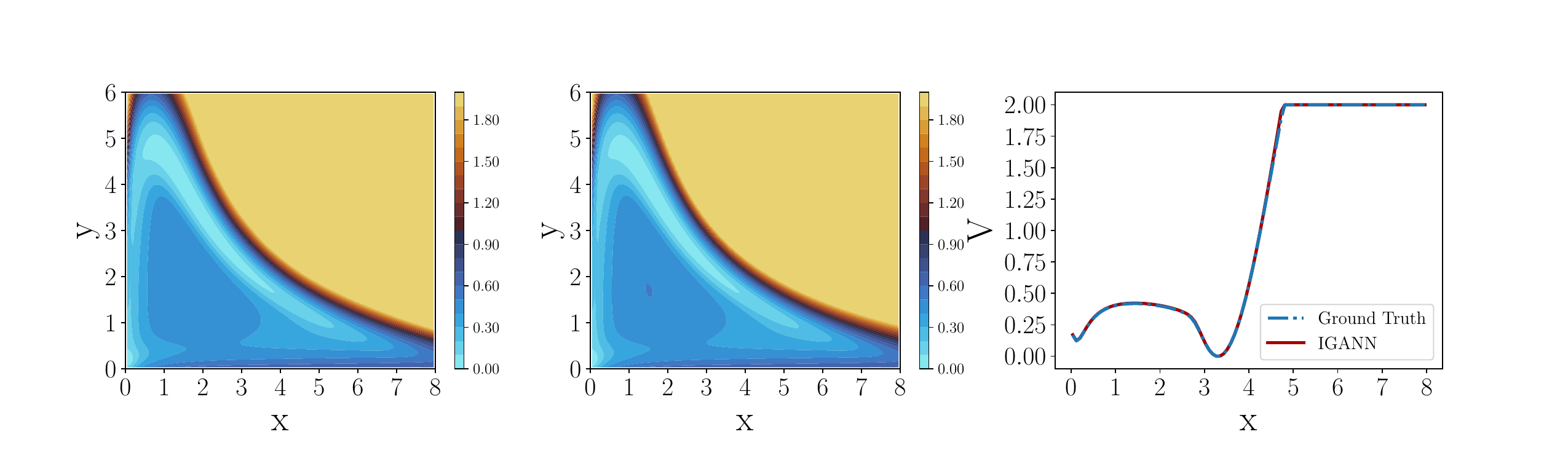}
	\caption{The learned potential $V_{\vtheta}$ compared with true potential for  $\varepsilon = 0.1$ (the IGANN prediction (left), exact potential (middle) and central line $y=2$ (right).}
	\label{fig: ex4.2.1 bio}
\end{figure}

\textbf{Biochemical oscillation network model.} For the limit cycle dynamics, the potential landscape will have a limit-cycle shape when the system temperature is low.
\begin{equation}
\left\{  
             \begin{array}{lr}  
             \mathrm{d} x = (\kappa(\frac{\alpha^2+x^2}{1+x^2} \frac{1}{1+y} - ax)) \mathrm{d} t + \sqrt{2\varepsilon} \mathrm{d}W^x, \\  
             \mathrm{d} y = (\frac{\kappa}{\tau_0}(b - \frac{y}{1+cx^2})) \mathrm{d}t + \sqrt{2\varepsilon} \mathrm{d}W^y, \\
             \end{array}  
\right.  
\end{equation}
where $\kappa=100$, $\alpha=a=b=0.1$, $c=100$, $\tau_0=5$, we conduct the computation for $\varepsilon=0.1$. In the first epoch, we use the LHS to sample 10000 samples in the larger domain $[-0.8, 12] \times [-0.8, 8]$, then we use these samples to simulate the SDE \eqref{eq:1} using the Euler-Maruyama scheme, and the resulting data serves as our training samples.  The learning potential is shown in Fig. \ref{fig: ex4.2.1 bio}, our IGANN method can effectively capture the features of the potential landscape, the rRMSE and rMAE are $5.63\%$ and $2.13\%$ respectively. These results are slightly better than the results in Ref \cite{lin2022computing}, where the rRMSE and rMAE were $8.97\% \pm 2.83\%$ and $6.63\% \pm 1.54\% $, respectively.
\FloatBarrier

\begin{figure}[H]
	\centering
	\includegraphics[width=1\linewidth]{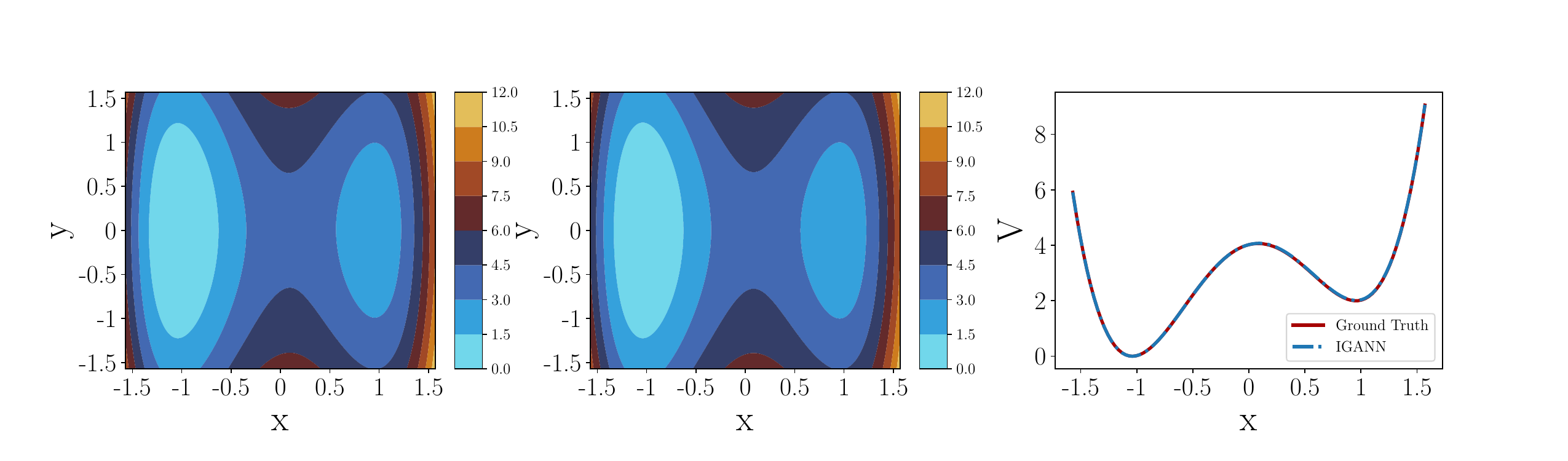}
	\caption{The learned potential $V_{\vtheta}$ compared with true potential for $\varepsilon = 0.04$ (the IGANN prediction (left), true potential (middle) and central line $y=0$ (right).}
	\label{fig: ex4.2.1 self-constructed}
\end{figure}

\textbf{Synthetic dynamical system example.} \label{ex:self-constructed} The dynamical driving force $\bm{f}$ is constructed in accordance with the given generalized potential $V=3(1-x^2)^2+y^2+x$ and the coefficient $\varepsilon$:
\begin{equation}
    \bm{f}(\vx) = \nabla \times \bm{f} - \frac{1}{\varepsilon}\nabla V(\vx) \times \bm{f} - \nabla V(\vx),
\end{equation}
where $F = (0, 0, 1)^{T}$. At each training epoch, we gather 2000 samples using LHS, and set $\beta=150$ with $\beta_{decay}=\num{3e-4}$. The learned potential $V_{\vtheta}$ is illustrated in Fig. \ref{fig: ex4.2.1 self-constructed}, it almost coincides with the true solution. In Table \ref{table: ex4.2.1 self-constructed}, we can see that it can learn the generalized potential well for different coefficients $\varepsilon$.

\FloatBarrier

\begin{figure}[H]
	\centering
	\includegraphics[width=1\linewidth]{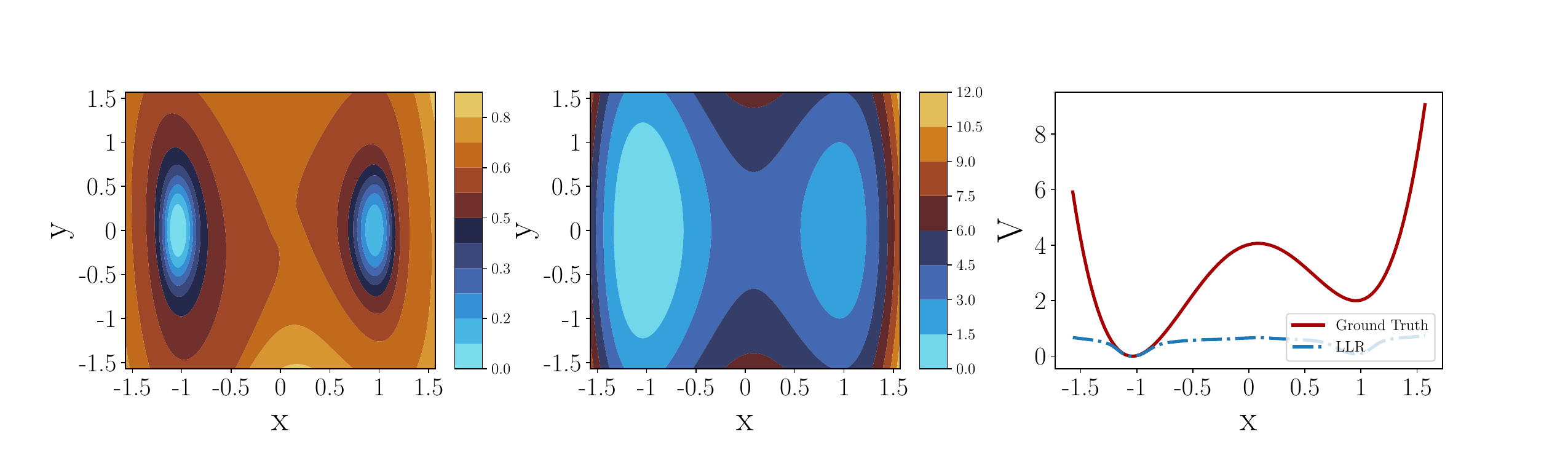}
	\caption{The learned potential $V_{\vtheta}$ using the LLR method with $\varepsilon = 0.1$ (the LLR prediction (left), true potential (middle) and central line $y=0$ (right).}
	\label{fig: ex4.2.1 Lin_fail}
\end{figure}
\FloatBarrier

We compare our IGANN method with the LLR method. For the synthetic dynamical system example, the solution computed using the LLR method \cite{lin2022computing} can not converge to the exact potential as illustrated in Fig. \ref{fig: ex4.2.1 Lin_fail}. The LLR method fails to recover the correct well.

\begin{table}[htbp]
\centering
\caption{Example \ref{sec: 4.2.1} self-constructed system}
\begin{tabular}{|p{2cm}<{\centering}|p{3cm}<{\centering}|p{3cm}<{\centering}|} \hline
&rRMSE (IGANN)&rMAE (IGANN)\\ \hline
$\varepsilon=0.1$ &  0.135\% & 0.130\% \\ \hline
$\varepsilon=0.07$&  0.611\%& 0.559\% \\ \hline
$\varepsilon=0.04$&  0.309\%& 0.285\%\\ \hline
\end{tabular}
\label{table: ex4.2.1 self-constructed}
\end{table}

\subsubsection{Three-dimensional systems} \label{sec 4.2.2}
\begin{figure}[H]
	\centering
	\includegraphics[width=1\linewidth]{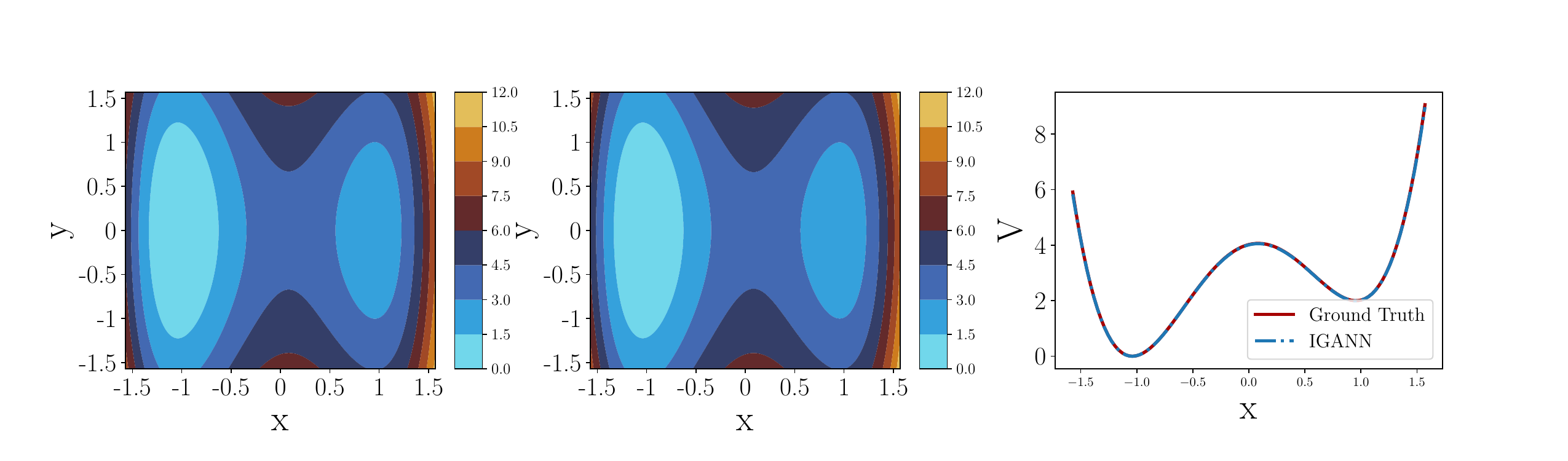}
	\caption{The learned potential $V_{\vtheta}$ compared with true potential for $\varepsilon = 0.04$ (the IGANN prediction (left), true potential (middle) and central line $y=0$ (right).}
	\label{fig: ex4.2.2 three dimensional}
\end{figure}

For a higher-dimensional non-gradient system, we consider a three-dimensional dynamical system with the dynamical driving force $\bm{f}$ constructed according to the given generalized potential $V=3(1-x^2)^2+y^2+z^2+x$ and the coefficient $\varepsilon$:
\begin{equation}
    \bm{f}(\vx) = \nabla \times \bm{f} - \frac{1}{\varepsilon}\nabla V(\vx) \times \bm{f} - \nabla V(\vx). \label{eq: f}
\end{equation}
During training, at each epoch, we gather 2000 samples using LHS. We set $\beta=40$ with $\beta_{decay}=\num{2e-4}$ for $\varepsilon=0.1$, $\beta=40$ with $\beta_{decay}=\num{1.5e-4}$ for $\varepsilon=0.07$, and $\beta=70$ with $\beta_{decay}=\num{1e-4}$ for $\varepsilon=0.04$. We demonstrate the cross-section of the learned potential $V_{\vtheta}$ in Fig. \ref{fig: ex4.2.2 three dimensional}, which almost coincides with the exact solution. The rRMSE and rMAE errors are provided in Table \ref{table: ex4.2.2}.
\FloatBarrier

\begin{table}[htbp]
\centering
\caption{Example \ref{sec 4.2.2}}
\begin{tabular}{|p{2cm}<{\centering}|p{3cm}<{\centering}|p{3cm}<{\centering}|} \hline
&rRMSE (IGANN)&rMAE (IGANN)\\ \hline
$\varepsilon=0.1$ &  1.79\% & 0.880\% \\ \hline
$\varepsilon=0.07$&  2.96\%& 0.969\% \\ \hline
$\varepsilon=0.04$&  2.03\%&  1.53\%\\ \hline
\end{tabular}
\label{table: ex4.2.2}
\end{table}

\subsubsection{Four-dimensional systems}
\begin{figure}[H]
	\centering
	\includegraphics[width=1\linewidth]{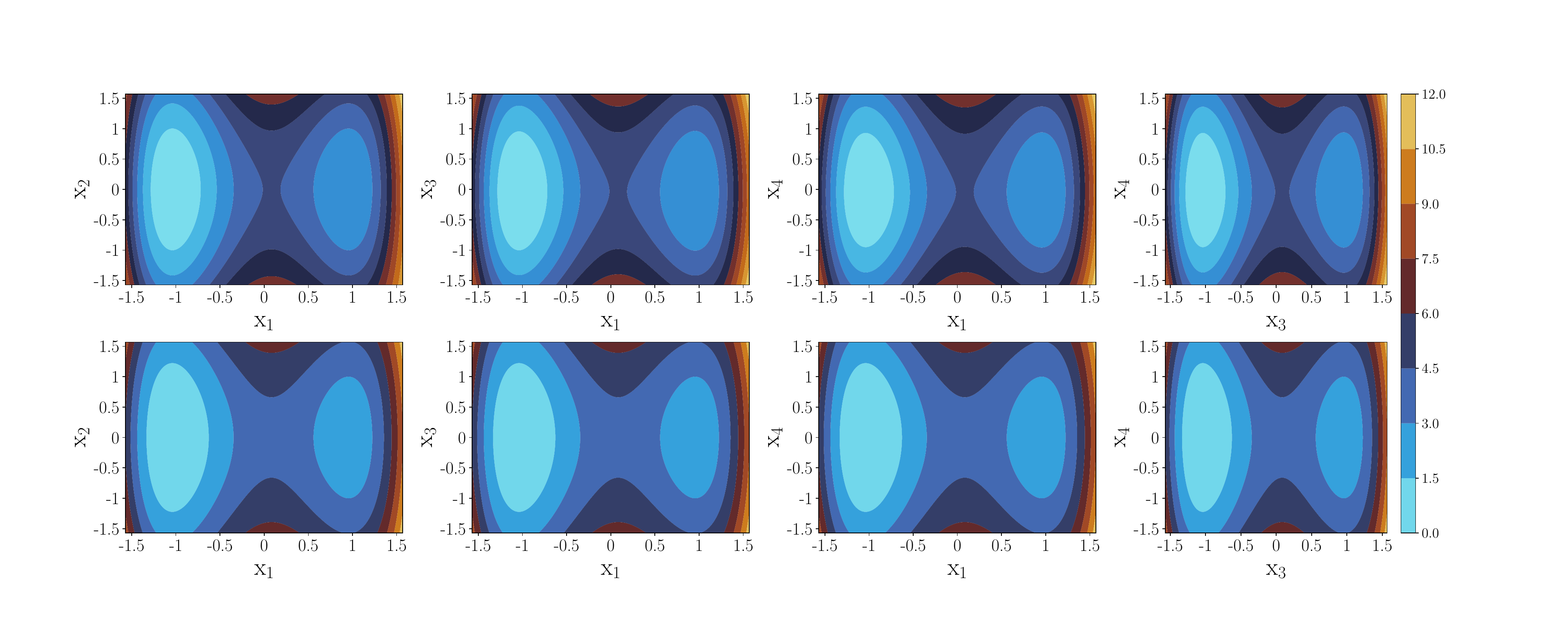}
	\caption{The learned potential $V_{\vtheta}$ compared with true potential for $\varepsilon = 0.04$ (the IGANN prediction (left), exact potential (middle) and central line $y=0$ (right).}
	\label{fig: ex4.2.3 1}
\end{figure}

\begin{figure}[H]
	\centering
	\includegraphics[width=1\linewidth]{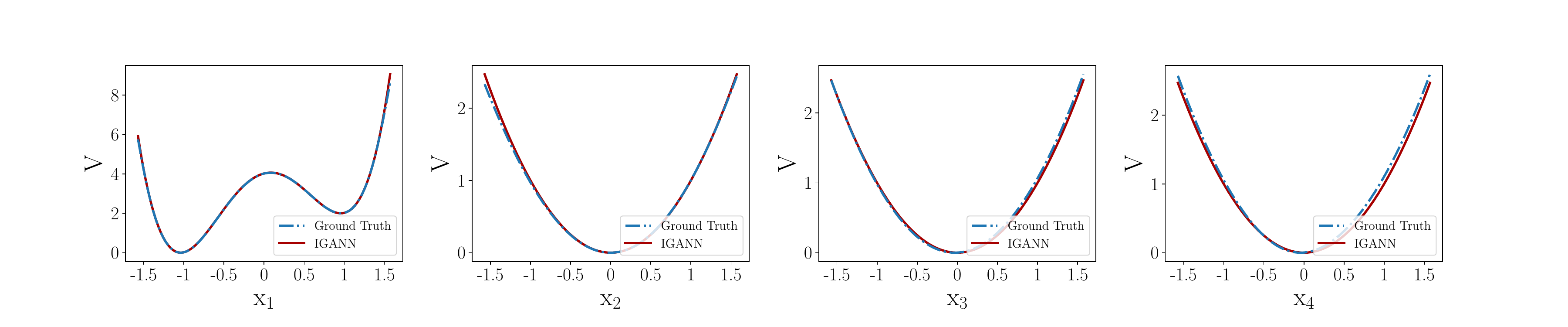}
	\caption{The learned potential $V_{\vtheta}$ compared with true potential for $\varepsilon = 0.04$ (the IGANN prediction (left), exact potential (middle) and central line $y=0$ (right).}
	\label{fig: ex4.2.3 2}
\end{figure}

In this section, we consider the four-dimensional dynamical system. The dynamical driving force $\bm{f}$ is constructed according to the Eq. \eqref{eq: f}, and the generalized potential is $V=3(1-x^2)^2+y^2+z^2+\omega^2+x$. We set $\beta=150$ with $\beta_{decay}=\num{2e-4}$ for $\varepsilon=0.04$. At each training epoch, we collect 2000 samples with LHS sampling. We show the cross section of the learned potential $V_{\vtheta}$ in Fig. \ref{fig: ex4.2.3 1} and Fig. \ref{fig: ex4.2.3 2}. As can be seen from the pictures, the solution almost coincides with the exact solution.
\FloatBarrier

\FloatBarrier

\subsubsection{Eight-dimensional systems}
\begin{figure}[H]
	\centering
	\includegraphics[width=1\linewidth]{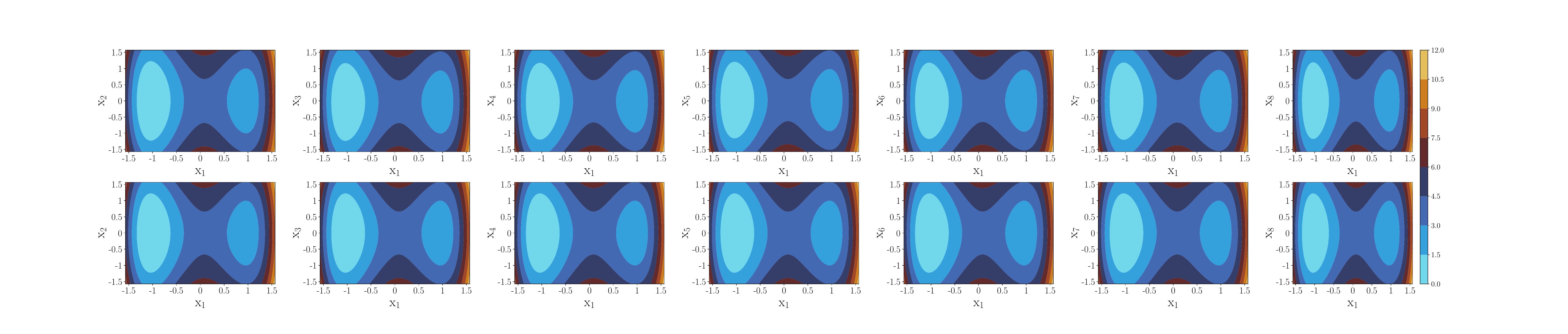}
	\caption{The upper pictures are learned potential $V_{\vtheta}$ (IGANN) and the lower pictures are exact potential $V$ for $\varepsilon = 0.04$.}
	\label{fig: ex4.2.4 suf}
\end{figure}

\begin{figure}[H]
	\centering
	\includegraphics[width=1\linewidth]{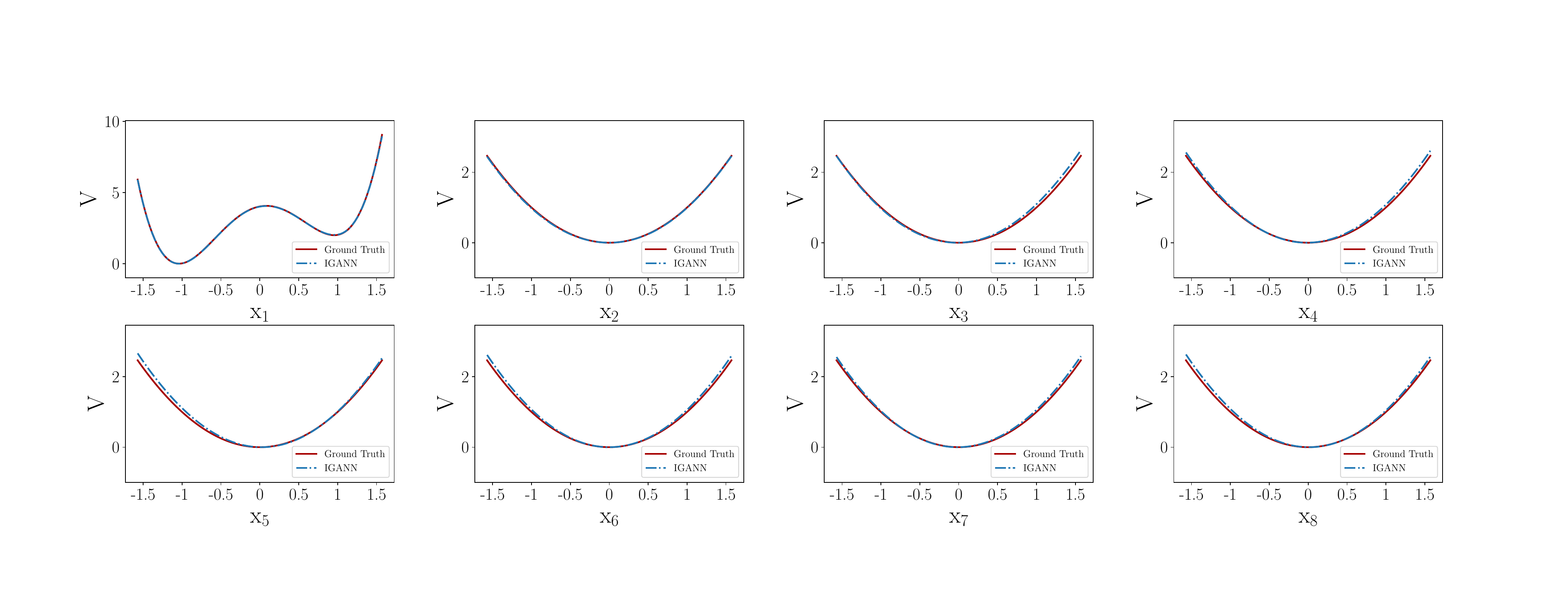}
	\caption{The learned potential $V_{\vtheta}$ (IGANN) compared with true potential for $\varepsilon = 0.04$.}
	\label{fig: ex4.2.4 line}
\end{figure}

In order to demonstrate the ability of our method for high dimension problems, an eight-dimensional dynamical system is constructed using Eq. \eqref{eq: f}, the generalized potential $V$ is given by $3(1-x_1^2)^2+\sum_{i=2}^{8}x_{i}^2+x$:
\begin{equation}
    \bm{f}(\vx) = \nabla \times \bm{f} - \frac{1}{\varepsilon}\nabla V(\vx) \times \bm{f} - \nabla V(\vx)
\end{equation}
We set $\beta=150$ with $\beta_{decay}=\num{2e-4}$ for $\varepsilon=0.04$. At each training epoch, we gather 2000 samples using LHS. We demonstrate the cross-section of the learned potential $V_{\vtheta}$ in Fig. \ref{fig: ex4.2.4 suf}, the upper picture is the learned solution in different coordinate planes $x_i- x_j$ by letting other coordinate components $x_k=0, k\neq i, j$ , and the lower picture is the exact solution in the corresponding plane. We can see the learned solution closely approaches the exact solution. For specific comparison, we present the solution on different axes $x_i$ in Fig. \ref{fig: ex4.2.4 line}, by letting the other coordinate components $x_j=0, j \neq i $, the well can be captured accurately.
\FloatBarrier

\FloatBarrier

\subsubsection{Ten-dimensional systems}\label{section:10d}
\begin{figure}[H]
	\centering
	\includegraphics[width=1\linewidth]{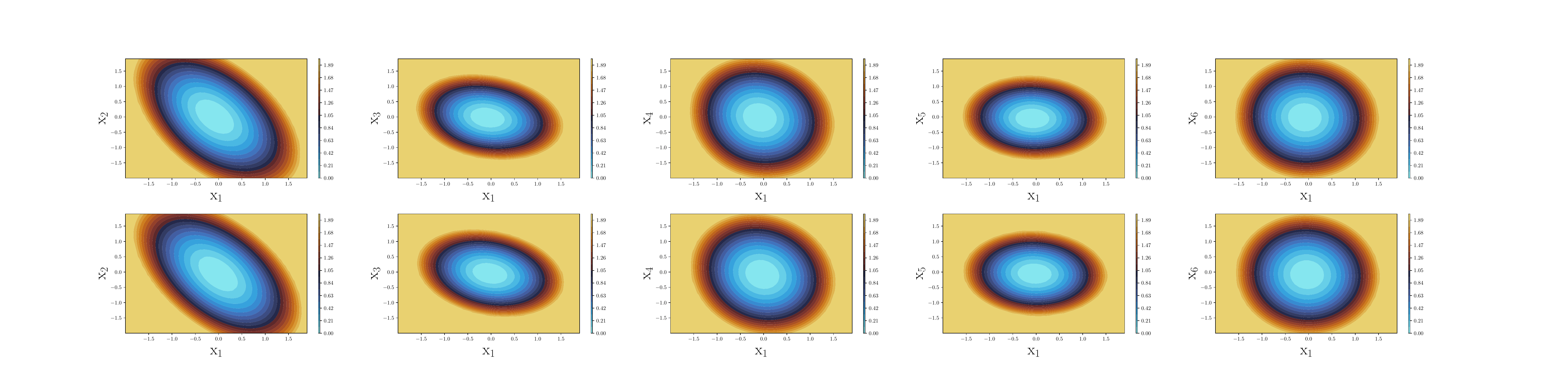}
	\caption{The learned potential $V_{\vtheta}$ (IGANN) compared with true potential for $\varepsilon = 0.1$.}
	\label{fig: 10d}
\end{figure}

We also test our method in a ten-dimensional system which is tested by LLR method \cite{lin2022computing}. The system is described by 

$$
\dot{\vx}=B \mathbf{h}\left(B^{-1} \vx\right)+\sqrt{2 \epsilon} B \mathbf{\xi}, \quad t>0,
$$
where $\mathbf{h}(\vy) = \left(h_1(\vy), \ldots, h_{10}(\vy)\right)^T$ is a vector field with
$$
\begin{aligned}
h_{2 k-1}(\vy) & =v_1\left(y_{2 k-1}, y_{2 k}\right):=-y_{2 k-1}+y_{2 k}\left(1+\sin y_{2 k-1}\right) \\
h_{2 k}(\vy) & =v_2\left(y_{2 k-1}, y_{2 k}\right):=-y_{2 k}-y_{2 k-1}\left(1+\sin y_{2 k-1}\right), \quad 1 \leq k \leq 5,
\end{aligned}
$$
$B=\left[b_{i, j}\right]$ is a $10 \times 10$ matrix given by
$$
b_{i, j}=\left\{\begin{array}{cl}
0.8, & \text { for } i=j=2 k-1,1 \leq k \leq 5 \\
1.25, & \text { for } i=j=2 k, 1 \leq k \leq 5 \\
-0.5, & \text { for } j=i+1,1 \leq i \leq 9 \\
0, & \text { otherwise }
\end{array}\right.
$$

and coefficient $\varepsilon=0.1$. We set $\beta=80$ with $\beta_{decay}=\num{3e-4}$. At each training epoch, we gather \num{1e4} samples using LHS. As shown in Fig. \ref{fig: 10d}, the results obtained by our IGANN method (upper) are completely consistent with the Ground Truth (lower).

\FloatBarrier
\section{Conclusion}
In this paper, we present the input gradient annealing neural network method, which is simple yet effective regardless of the temperature of the dynamical system. This method can help compute the exact solution by adding a negative penalty  term in the loss function. We examine its effectiveness in various examples, including gradient systems and non-gradient systems across different temperatures. In all examples, numerical solutions agree well with the reference solutions. Furthermore, this method does not require any prior knowledge of the free energy landscape or the force field, enabling us to study low-temperature systems.

\bibliographystyle{plain}
\bibliography{sample}

\end{document}